\documentclass{bcp}

 \usepackage{amsmath,amsthm}
\usepackage{amssymb,latexsym}

\newtheorem{theorem}{Theorem}[section]

\newtheorem{question}[theorem]{Question}
\newtheorem{conjecture}[theorem]{Conjecture}
\newtheorem{remark}[theorem]{Remark}

\theoremstyle{definition}

\newcommand{\what}[1]{\widehat{#1}}
\newcommand{\comp}{\raisebox{.2ex}{${\scriptstyle\circ}$}}
\newcommand{\con}{\ast}
\newcommand{\cross}{\!\times\!}
\newcommand{\inprod}[2]{\left\langle #1 \mid #2 \right\rangle}
\newcommand{\mult}{\cdot}

\newcommand{\id}{\operatorname{id}}

\newcommand{\supp}{\operatorname{supp}}

\newcommand{\norm}[1]{\left\Vert#1\right\Vert}

\newcommand{\fA}{\mathcal{A}}
\newcommand{\fB}{\mathcal{B}}
\newcommand{\fC}{\mathcal{C}}
\newcommand{\fH}{\mathcal{H}}
\newcommand{\fI}{\mathcal{I}}
\newcommand{\fM}{\mathcal{M}}
\newcommand{\fN}{\mathcal{N}}
\newcommand{\fV}{\mathcal{V}}
\newcommand{\fW}{\mathcal{W}}
\newcommand{\fX}{\mathcal{X}}
\newcommand{\fY}{\mathcal{Y}}
\newcommand{\fZ}{\mathcal{Z}}

\newcommand{\Cee}{\mathbb{C}}
\newcommand{\Ree}{\mathbb{R}}
\newcommand{\En}{\mathbb{N}}

\newcommand{\bloneg}{\ensuremath{\mathrm{L}^1(G)}}
\newcommand{\blone}[1]{\ensuremath{\mathrm{L}^1(#1)}}
\newcommand{\bltwog}{\ensuremath{\mathrm{L}^2(G)}}
\newcommand{\bltwo}[1]{\ensuremath{\mathrm{L}^2(#1)}}
\newcommand{\cstarg}{\ensuremath{\mathrm{C}^*(G)}}
\newcommand{\falg}{\ensuremath{\mathrm{A}(G)}}
\newcommand{\falh}{\ensuremath{\mathrm{A}(H)}}
\newcommand{\fal}[1]{\ensuremath{\mathrm{A}(#1)}}
\newcommand{\falgk}{\ensuremath{\mathrm{A}(G/K)}}
\newcommand{\ffalg}{\ensuremath{\mathrm{A}_f(G)}}
\newcommand{\fgamg}{\ensuremath{\mathrm{A}_\gamma(G)}}
\newcommand{\fdelg}{\ensuremath{\mathrm{A}_\Delta(G)}}
\newcommand{\fcbal}[1]{\ensuremath{\mathrm{A}_{cb}(#1)}}
\newcommand{\fsalg}{\ensuremath{\mathrm{B}(G)}}
\newcommand{\fsalh}{\ensuremath{\mathrm{B}(H)}}
\newcommand{\fsal}[1]{\ensuremath{\mathrm{B}(#1)}}
\newcommand{\fsralg}{\ensuremath{\mathrm{B}_r(G)}}
\newcommand{\fsoalg}{\ensuremath{\mathrm{B}_0(G)}}
\newcommand{\fsoal}[1]{\ensuremath{\mathrm{B}_0(#1)}}
\newcommand{\ideal}{\mathrm{I}}
\newcommand{\matm}[1]{\ensuremath{M}_m(#1)}
\newcommand{\matn}[1]{\ensuremath{M}_n(#1)}
\newcommand{\measg}{\ensuremath{\mathrm{M}(G)}}
\newcommand{\meas}[1]{\ensuremath{\mathrm{M}(#1)}}
\newcommand{\multcb}{\ensuremath{\mathrm{M}_{cb}}}
\newcommand{\vng}{\ensuremath{\mathrm{VN}(G)}}
\newcommand{\vnh}{\ensuremath{\mathrm{VN}(H)}}

\begin{document}

\mathclass{Primary 43-02, 43A30, 46H25, 46L07; 
Secondary 43A07, 43A20, 43A10, 43A77, 22D10, 46J20, 43A85.}
\abbrevauthors{N. Spronk}
\abbrevtitle{Amenability properties}

\title{
Amenability properties of Fourier algebras \\
and Fourier-Stieltjes algebras:  a survey
}


\author{Nico Spronk}
\address{Department of Pure Mathematics, University of Waterloo \\
200 University Ave.\ W., Waterloo, Ontario, N2L\,3G1, Canada \\
E-mail: {\tt nspronk@uwaterloo.ca}}

\maketitlebcp

\begin{abstract}
Let $G$ be a locally compact group, and let \falg\ and \fsalg\ denote its
Fourier and Fourier-Stieltjes algebras.  These algebras are dual objects
of the group and measure algebras, \bloneg\ and \measg, in a sense which 
generalizes the Pontryagin duality theorem on abelian groups.  We
wish to consider the amenability properties of \falg\ and \fsalg\ and compare
them to such properties for \bloneg\ and \measg.  For us, ``amenability
properties'' refers to amenability, weak amenability, and biflatness,
as well as some  
properties which are more suited to special settings, such as the 
hyper-Tauberian property for semisimple commutative Banach algebras.
We wish to emphasize that the theory of operator spaces and completely
bounded maps plays an indispensable role when studying \falg\ and
\fsalg.  We also show some applications of amenability theory to problems
of complemented ideals and homomorphisms.

\end{abstract}

\section{Introduction} 

This article, as the title literally suggests, is a survey of the
amenability results around Fourier and Fourier-Stieltjes algebras which are
known to the author.
These are displayed in comparison to the results on group and measure
algebras.  In particular, I intend to highlight the indispensable role
which operator spaces play in this theory.

The scope of this article has been purposely restricted to only 
amenability properties of Fourier and 
Fourier-Stieltjes algebras, and to the motivating results in their dual objects,
group and measure algebras.  
There are, thus, clear omissions of
discussions of Herz--Fig\`{a}-Talamanca algebras $\mathrm{A}_p(G)$,
in general, and their amenability properties; and of locally
compact quantum groups or even Kac algebras, though they do provide
a convenient language for discussing the duality.  However, 
limited liberty is taken to indicate some relevant literature on such topics,
as they relate to the present article.  
The scope of this article has been restricted mainly because the topic
has reached a certain maturity, though there are still interesting
open problems on Fourier and Fourier-Stieltjes algebras.  
Moreover,  it is hoped that this limitation in scope
improves coherence and depth, though at the cost of breadth.

We acknowledge that this is far from being the first survey of amenability
problems in Fourier and Fourier-Stieltjes algebras.  See surveys of
Runde~\cite{rundeS,rundeS1} and Kaniuth and Lau~\cite{kaniuthlS}.  However,
the subject has advanced significantly, even since the latter of these
surveys.  Moreover, our focus differs from the foci of these articles, thus we
offer a complement to them.

\subsection{Some Banach algebras of harmonic analysis}
Let $G$ be a locally compact group.  The {\it group} and {\it measure} 
convolution algebras of $G$, are denoted by
\bloneg\ and \measg, respectively.  These classical objects
are defined in many texts; see the treatise \cite{hewittrI} for definitions
and historical details.  These algebras represent complete invariants
for the underlying group in the sense that if $H$ is another locally compact
group for which $\bloneg\cong\blone{H}$, or $\measg\cong\meas{H}$, 
isometrically isomorphically, 
then $G\cong H$ isomorphically and homeomorphically, 
as shown by Wendel~\cite{wendel}.
We may, and will, regard \bloneg\ as the ideal in \measg\ of measures
which are absolutely continuous with respect to the Haar measure.

If $G$ is abelian, it admits a dual group
$\what{G}$, and the Fourier-Stieltjes transform $\mu\mapsto\hat{\mu}:
\meas{\what{G}}\to\fC_b(G)$ (bounded continuous functions on $G$)
restricts to the Gelfand transform $f\mapsto\hat{f}:\blone{\what{G}}
\to\fC_0(G)$ (continuous functions vanishing at infinity).  Their respective
images \fsalg\ and \falg, each endowed with the norm which makes the transform 
an isometry, are thus Banach algebras of functions, which are invariants
for $G$.  It follows from Jordan decomposition in \meas{\what{G}}, 
and Bochner's theorem, that \fsalg\ is the span of 
continuous positive definite functions; and then a version of the
Gelfand-Naimark-Segal construction shows that \fsalg\ consists of the family
of continuous matrix coefficients of all unitary representations on $G$.
Since each element of \bloneg\ factors as a pointwise product of two
square-integrable functions, i.e.\ $\blone{\what{G}}=\bltwo{\what{G}}\mult
\bltwo{\what{G}}$, it follows from the Plancherel theorem that
$\falg=\bltwog\con\bltwog$.

Eymard~\cite{eymard} showed, using representation theory, how to construct
generalizations of \falg\ and \fsalg\ for any locally compact group.
The {\it Fourier-Stieltjes algebra} \fsalg\ is the space of continuous
matrix coefficients of unitary representations of $G$: 
$s\mapsto\inprod{\pi(s)\xi}{\eta}$ for any weak operator continuous
unitary representation $\pi:G\to\fB(\fH)$, $\xi,\eta\in\fH$.
The {\it Fourier algebra} \falg\ is the space of matrix coefficients of the
left regular representation: $s\mapsto\inprod{\lambda(s) f}{g}
=\bar{g}\con\check{f}(s)$, $f,g\in\bltwog$, where 
$\lambda:G\to\fB(\bltwog)$ is given by $\lambda(s)f(t)=f(s^{-1}t)$,
and $\check{f}(t)=f(t^{-1})$.  \fsalg\ may be identified as the dual space
of the universal C*-algebra $\cstarg$, of $G$, and this norm makes
\fsalg\ a Banach algebra, under pointwise operations. $\falg$
is a closed ideal of \fsalg\, a fact which may be verified by showing
that $\falg$ is generated by the compactly supported elements of \fsalg,
or, alternatively, by an application of Fell's ``absorption principle''
\cite{fell}.
These are semisimple commutative Banach algebras, and $G$ is the Gelfand
spectrum of \falg.  Walter~\cite{walter} showed the analogue of Wendel's 
theorem: if $H$ is another locally compact group then $\falg\cong\falh$,
or $\fsalg\cong\fsalh$, isometrically isomorphically, exactly
when $G\cong H$ isomorphically and homeomorphically.
Walter's result foreshadows the use of operator space methods, in
that he employs Kadison's characterizations of surjective isometries
between C*-algebras \cite{kadison}.

It will be useful to identify various subalgebras of \fsalg.
We let \fsralg\ denote the weak*-closure of \falg\ in \fsalg.
It is a result of Hulanicki~\cite{hulanicki} that
\fsralg=\fsalg\ if and only if $G$ is amenable.  Let \ffalg\
denote the norm-closed subspace of \fsalg\ generated by matrix
coefficients of continuous finite dimensional representations.
Then $\ffalg\cong\fal{G^{ap}}$ where $G^{ap}$ is the almost periodic
compactification of $G$.  We note that if $G$ is abelian, $\ffalg$
is the image under the Fourier transform of the closed span of 
Dirac measures $\ell^1(\what{G})$
in $\meas{\what{G}}$.  We let the Rajchman algebra be given by
$\fsoalg=\fsalg\cap\fC_0(G)$.  For any non-compact abelian group, it is known that
$\fsoalg\supsetneq\falg$; see \S 6 of \cite{grahamm}.

\subsection{Amenability and related properties}\label{ssec:amenability}
If $\fA$ is a Banach algebra, a Banach space $\fX$ is called
a (contractive) Banach $\fA$-bimodule if there are bounded
(contractive) maps, a homomorphism $\fA\to\fB(\fX):a\mapsto(x\mapsto a\mult x)$
and an anti-homomorphism 
$\fA\to\fB(\fX):a\mapsto(x\mapsto x\mult a)$, with commuting ranges.
The adjoints of these maps make the dual space $\fX^*$ into
a dual Banach $\fA$-bimodule.  A linear map $D:\fA\to\fX$
is called a {\it derivation} if $D(ab)=a\mult D(b)+D(a)\mult b$
for $a,b$ in $\fA$.  Inner derivations are those of the form 
$D(a)=a\mult x-x\mult a$ for some $x$ in $\fX$.

In his seminal memoir \cite{johnsonM},
Johnson defined $\fA$ to be {\it amenable} if, for every
dual Banach $\fA$-bimodule $\fX^*$, every bounded derivation
$D:\fA\to\fX^*$ is inner.  He further showed in \cite{johnson72}
that amenability is equivalent to a certain averaging property:
$\fA$ is amenable if and only if it admits a {\it bounded approximate
diagonal}, i.e.\ a bounded net $(d_\alpha)$ in $\fA\otimes^\gamma\fA$
(projective tensor product) such that
\[
m(d_\alpha)a,am(d_\alpha)\to a\quad\text{ and }\quad
a\mult d_\alpha-d_\alpha\mult a\to 0
\]
where $m:\fA\otimes^\gamma\fA\to\fA$ is the multiplication map
$m(a\otimes b)=ab$, and $a\mult(b\otimes c)=(ab\times c),
(b\otimes c)\mult a=b\otimes(ca)$.  Bounded approximate
diagonals allow for a quantitative measurement of amenability:  we say
$\fA$ is $C$-amenable if it admits a bounded approximate diagonal
of norm at most $C$, and let the {\it amenability constant} of
$\fA$, $C_\fA$ denote the infimum of such $C$.
Bounded approximate diagonals
are also useful in the following application of amenability
due to Helemski\u{\i} (see the monograph \cite{helemskii} or \cite{curtisl};
the role of the bounded approximate diagonal is apparent in 
\cite{rundeB}):  if $\fA$ is amenable and $\fX$ is a Banach $\fA$-bimodule
with a boundedly complemented subspace $\fY$, which is also an
$\fA$-bimodule, then there is a bounded projection $P:\fX^*\to\fY^\perp$
which is an $\fA$-bimodule map.

We say $\fA$ is {\it contractible} or {\it super-amenable} if
all bounded derivations into (not necessarily dual) $\fA$-bimodules
are inner.  A contractible C*-algebra is necessarily 
finite dimensional (see \cite{rundeB}, for example).

We say $\fA$ is a {\it dual} Banach algebra if there is a closed
$\fA$-bimodule $\fA_*\subset\fA^{**}$ such that $(\fA_*)^*\cong\fA$.
If every weak*-weak* continuous derivation $D:\fA\to\fX^*$ is inner,
we say $\fA$ is {\it Connes amenable}.  This terminology first appeared in
\cite{runde01}, and honours Connes's version of amenability
for von Neumann algebras.

A Banach algebra is {\it weakly amenable} if every bounded derivation
$D:\fA\to\fA^*$ is inner.  For the case of commutative $\fA$, this terminology
was coined in \cite{badecd}, where it was shown that any derivation
$D:\fA\to\fX$ for symmetric $\fX$ (right and left actions coincide)
vanishes, exactly when $\fA$ is weakly amenable.  Given a character
character $\chi:\fA\to\Cee$, a non-zero linear functional $d\in\fA^*$ such that 
$d(ab)=\chi(a)d(b)+d(a)\chi(b)$ is called a point derivation.  If a non-zero
derivation exists, the derivation $D:\fA\to\fA^*$, $D(a)=d(a)\chi$,
shows that $\fA$ cannot be weakly amenable.

We indicate some related properties due to Helemski\u{\i}~\cite{helemskii}.
$\fA$ is {\it biprojective} if there is a bounded $\fA$-bimodule
map $S:\fA\to\fA\otimes^\gamma\fA$
such that $m\comp S=\id_\fA$.  $\fA$ is {\it biflat} if there is a 
bounded $\fA$-bimodule map $T:(\fA\otimes^\gamma\fA)^*\to\fA^*$
such that $m^*\comp T=\id_{\fA^*}$.  A biflat algebra is automatically
weakly amenable.  We also have that $\fA$ is amenable if and only if
it is biflat and has a bounded approximate identity.

\subsection{Commutative Banach algebras} \label{sec:commutative}
We let, for this section, $\fA$ be a regular function algebra on a locally compact
space $X$, i.e.\ $\fA$ is a semisimple commutative Banach algebra for which
there is a contractive inclusion $\fA\subset\fC_0(X)$
whose image contains functions which separate compact sets from
disjoint closed sets.  If $E\subset X$ is closed we let
\begin{gather*}
\ideal_\fA(E)=\{u\in\fA:u|_E=0\},\quad
\ideal_\fA^c(E)=\{u\in\ideal_\fA:\supp u\text{ is compact}\} \\
\text{and }\ideal_\fA^0(E)=\{u\in\ideal_\fA^c:\supp u\cap E=\varnothing\}.
\end{gather*}
We note that $\ideal_\fA^0(E)$ is the smallest ideal of elements of $\fA$
which vanish on $E$ (see \cite{reiter}), and $\ideal_\fA(E)$ is the largest.
We say that $E$ is a {\it spectral} set for $\fA$ if 
$\overline{\ideal_\fA^0(E)}=\ideal_\fA(E)$,
{\it approximable} if $\ideal_\fA(E)$ has a bounded approximate identity,
{\it locally spectral} if $\overline{\ideal_\fA^0(E)}\supset\ideal_\fA^c(E)$, and
{\it essential} if $\overline{\ideal_\fA(E)^2}=\ideal_\fA(E)$. 
The Tauberian property for $\fA$ is that $\varnothing$ is spectral for $\fA$.

There is a notion of support of a functional $\mu$ in $\fA^*$, which
generalizes that for a measure.  Using this, Samei \cite{samei06}
devised an amenability-type property for some Banach function algebras.
$\fA$ is called {\it hyper-Tauberian} if every bounded linear local map
$T:\fA\to\fA^*$, i.e\ $\supp(Tu)\subset\supp u$ for $u$ in $\fA$, is
an $\fA$-module map.  He proved that the hyper-Tauberian condition
implies weak amenability.

We have the following relationship between amenability conditions and spectral
conditions.

\begin{theorem}\label{theo:spectral}
Let $\fA$ be a regular function algebra on $X$ for
which $\fA\otimes^\gamma\fA$ is semi-simple.  Then
$\fA\otimes^\gamma\fA$ is regular on $X\cross X$ and
$\ideal_{\fA\otimes^\gamma\fA}(\Delta)=\ker m$ where
$\Delta=\{(x,x):x\in X\}$.  Moreover we have
\begin{itemize}
\item[\rm(i)] $\fA$ is amenable if and only if
$\varnothing$ is approximable for $\fA$ and $\Delta$ is approximable
for $\fA\otimes^\gamma\fA$;
\item[\rm(ii)] $\fA$ is hyper-Tauberian if and only if
$\varnothing$ is spectral for $\fA$ and $\Delta$ is locally spectral
for $\fA\otimes^\gamma\fA$; 
\item[\rm(iii)] if $\varnothing$ is approximable for $\fA$, then
$\fA$ is weakly amenable if and only
$\Delta$ is essential for $\fA\otimes^\gamma\fA$;
\item[\rm(iii')] $\fA$ is weakly amenable if $\varnothing$ is essential 
for $\fA$ and $\Delta$ is spectral for $\fA\otimes^\gamma\fA$.
\end{itemize}
\end{theorem}

We note that (i) is a specialized version of a famous splitting result
of \cite{helemskii}; see also \cite{curtisl}.  The result
(ii) is a mild generalization
of a result in \cite{samei06}; there, the only extra assumption is that $X$ is the
spectrum of $\fA$.  The result (iii) is in \cite{groenbaek}; (iii') is an easy corollary
of the same theorem of \cite{groenbaek} in which (iii) is stated.

\begin{question}\label{ques:hypertaub}
\begin{itemize}
\item[\rm (i)] Is it the case that amenability implies hyper-Tauberianness?
\item[\rm (ii)] Is there an example of an approximable set which
is not (locally) spectral?
\end{itemize}
\end{question}

A positive answer to (ii) would dismiss an obvious route to proving (i).

\section{Amenability properties of Banach algebras of harmonic analysis}

\subsection{Group and measure algebras} \label{sec:blonemeas}
Let $G$ be a locally compact group.  The property of
{\it amenability} for groups is well known; see the monograph of
Paterson~\cite{paterson}.

\begin{theorem}\label{theo:johhel}
The following are equivalent:
\begin{itemize}
\item[\rm (i)] $G$ is an amenable group;
\item[\rm (ii)] \bloneg\ is amenable; and
\item[\rm (iii)] \bloneg\ is biflat.
\end{itemize}
\end{theorem}

The equivalence of (i) and (ii) is a famous and motivating theorem
of Johnson~\cite{johnsonM}.  The equivalence of (i) and (iii) is due
to Helemski\u{\i}~\cite{helemskii}.  Note that since $\bloneg$ always has a bounded 
approximate identity, (ii) and (iii) are equivalent.  It is well-known
that \bloneg, which injects densely into a C*-algebra (say $\cstarg$), is
contractible if and only if $G$ is finite.

\begin{theorem}\label{theo:bloneweakam}
$\bloneg$ is always weakly amenable.
\end{theorem}

This result is due to Johnson~\cite{johnson91}; a simpler proof may be
found in \cite{despicg}.  Helemski\u{\i}~\cite{helemskii} proved the biprojectivity
result, below.

\begin{theorem}\label{theo:blonebiproj}
$\bloneg$ is biprojective if and only if $G$ is compact.
\end{theorem}

The characterization of amenability for measure algebras
is due to Dales, Ghahramani and Helemski\u{\i}~\cite{dalesgh}.

\begin{theorem}\label{theo:daghhe}
The following are equivalent:
\begin{itemize}
\item[\rm (i)] $\measg$ is weakly amenable;
\item[\rm (ii)] $\measg$ admits no point derivations; and
\item[\rm (iii)]  $G$ is discrete (and hence $\measg=\ell^1(G)$).
\end{itemize}
\end{theorem}

It is an immediate consequence that $\measg$ is amenable if and
only if $G$ is discrete and amenable.  For abelian groups, Theorem~\ref{theo:daghhe}
was established by Brown and Moran~\cite{brownm}.
We note that since $\measg$ is unital, if it is biprojective it is then contractible. 
Moreover since \measg\ injects densely into a C*-algebra
(for example the closure of $\lambda(\measg)$ in \vng), \measg\
is contractible if and only if $G$ is finite.

We close this section with a result of Runde~\cite{runde03}.

\begin{theorem}\label{theo:connes}
\measg\ is Connes amenable if and only if $G$ is amenable.
\end{theorem}

\subsection{Fourier and Fourier-Stieltjes algebras} \label{sec:falgfsalba}
For an abelian locally compact
group $G$, we have that $G$ is compact if and only if $\widehat{G}$ is discrete;
and, by Pontryagin duality, $G$ is discrete if and only if $\what{G}$ is compact.
Given this, it was long expected for general locally compact $G$, the results
for \falg\ and \fsalg\ would parallel those for \bloneg\ and \measg.
In particular, Leptin~\cite{leptin} showed that $\falg$ has a bounded
approximate identity if and only if $G$ is amenable. In particular,
if \falg\ is amenable, then we must have $G$ is amenable.  It was, for
a long time, thought that the converse must hold.
Hence the following result of Johnson~\cite{johnson94} was a surprise.

\begin{theorem}\label{theo:falgnwa}
For the compact group $G=\mathrm{SO}(3)$, \falg\
is not weakly amenable.
\end{theorem}

This was proved by taking a ``twisted'' convolution map
$u\otimes v\mapsto u\con\check{v}:\falg\otimes^\gamma\falg\to\falg$, 
and identifying its range
$\fgamg$ as a Banach algebra in its own right.  A non-zero point derivation
was found on $\fgamg$, which was used to show that
$\Delta$ is not essential for $\falg\otimes^\gamma\falg\to\falg$
(see Theorem~\ref{theo:spectral}).  Plymen~\cite{plymen}
showed that such a point derivation can be found for
$\fgamg$ of any compact semi-simple Lie group.  

The actual characterization of amenability is due to Forrest and 
Runde~\cite{forrestr} (see also a more qualitative version in \cite{runde06}).
We say $G$ is virtually abelian if it admits an abelian subgroup of finite index.

\begin{theorem}\label{theo:runfor}
The following are equivalent:
\begin{itemize}
\item[\rm (i)] \falg\ is amenable;
\item[\rm (ii)] $\check{\Delta}=\{(s,s^{-1}):s\in G\}\in\Omega(G)$, where
$\Omega(G)$ denotes the smallest ring of subsets containing all cosets; and
\item[\rm (iii)] $G$ is virtually abelian. 
\end{itemize}
Moreover, \fsalg\ is amenable if and only if $G$ is compact and virtually
abelian.
\end{theorem}

It is interesting to note that operator space techniques were used
in the proof (and that this proof is intimately related to Theorem~\ref{theo:iliespr},
below).  It seems that (ii)$\Leftrightarrow$(iii) can be proved by purely group theoretic
techniques, though the author knows of no reference for this.  
The implication (iii)$\Rightarrow$(i) can be found in \cite{laulw}.
We note that if \fsralg\ is amenable, it must have a bounded approximate identity
whose cluster point must be the constant function $1$; this means that 
\fsralg=\fsalg.  Contractibility of \falg\ is equivalent to finiteness of $G$;
indeed \falg\ injects densely into the commutative C*-algebra $\fC_0(G)$.

The characterization for weak amenability is not yet entirely known, but we have
strong partial results.

\begin{theorem}\label{theo:falwa}
\begin{itemize}
\item[\rm (i)] If the connected component of the identity $G_e$ of
$G$ is abelian, then \falg\ is hyper-Tauberian (hence weakly amenable).
\item[\rm (ii)] If $G$ contains a non-abelian connected compact subgroup,
then \falg\ is not weakly amenable.
\item[\rm (iii)] If $G$ is compact and connected, \falg\ is hyper-Tauberian
if and only if $G$ is abelian.
\end{itemize}
\end{theorem}

We note that connected [SIN] (small invariant neighbourhood) and
[MAP] (maximally almost periodic) groups (see tables in \cite{palmerB})
contain non-abelian connected compact subgroups whenever they, themselves,
are non-abelian.
Forrest and Runde~\cite{forrestr} proved (i) for weak amenability, while
Samei~\cite{samei06} proved it for the hyper-Tauberian property.
Results (ii) and (iii) were proved in \cite{forrestss2}; they follow from the 
next result in the same article which we give below.  
See Section \ref{sec:commutative} for notation
and terminology.

\begin{theorem}\label{theo:falwac}
Let $G$ be a compact group.  Then the following are equivalent:
\begin{itemize}
\item[\rm (i)] $G_e$ is abelian;
\item[\rm (ii)] \falg\ is hyper-Tauberian;
{\rm (ii')} $\Delta$ is spectral for $\falg\otimes^\gamma\falg$;
\item[\rm (iii)] \falg\ is weakly amenable; and
{\rm (iii')} $\Delta$ is essential $\falg\otimes^\gamma\falg$.
\end{itemize}
\end{theorem}

The proof of this theorem  builds on that of Theorem
\ref{theo:falgnwa} and the associated result of Plymen, and uses
the characterization of connected compact groups from \cite{price}.
The equivalences (ii)$\Leftrightarrow$(ii') and
(iii)$\Leftrightarrow$(iii') are direct applications of Theorem \ref{theo:spectral}.

\begin{question}
Is it the case that \falg\ is weakly amenable only if $G_e$ is abelian?
Is the Fourier algebra weakly amenable for any of the groups
$\mathrm{SL}_2(\Ree)$, the $ax+b$-group, or any of the Heisenberg groups?
\end{question}

Given Theorems \ref{theo:daghhe}, \ref{theo:runfor} and \ref{theo:falwa}, 
it is reasonable
to expect that \fsalg\ ought to be weakly amenable if and only if $G$ is compact
and $G_e$ is abelian; compare this with Theorem~\ref{theo:fell}~(iii),
below.

The best results known for biflatness and biprojectivity are recently
due to Runde~\cite{runde09}; though some components of this
theorem were developed by Aristov~\cite{aristov}.

\begin{theorem}\label{theo:falgbiflat}
\begin{itemize}
\item[\rm (i)] If \falg\ is biflat then either {\rm (a)} $G$ is virtually abelian, or
{\rm (b)} $G$ is non-amenable and does not contain a discrete copy of 
the free group on two generators.
\item[\rm (ii)] If \falg\ is biprojective then $G$ is discrete and
one of {\rm (a)}, or {\rm (b)}, in {\rm (i)} above, holds.
\end{itemize}
Conversely, if {\rm (a)} holds (and $G$ is discrete) then \falg\ is biflat (biprojective).
\end{theorem}

It seems unlikely that \falg\ could be biflat for any non-amenable group,
but little is known about harmonic analysis on 
non-amenable groups not containing free groups;
see Ol'shanskii and Sapir \cite{olshanskiis} for an up-to-date treatment
of the groups themselves.  Much as is the case for \measg, above,
\fsalg\ is biprojective if and only if $G$ is finite:  \fsalg\ has a commutative C*-algebra,
the Eberlein algebra, as its uniform closure in $\fC_b(G)$.

Results on Connes amenability for \fsalg\ are also incomplete.  The best
results published are due to Runde~\cite{runde04}.

\begin{theorem}\label{theo:fsalgconnes}
\begin{itemize}
\item[\rm (i)] If $G$ is either discrete and amenable, or is a product
of finite groups, then \fsalg\ is Connes amenable if and only if 
$G$ is virtually abelian.
\item[\rm (ii)] If $G$ is discrete, the \fsralg\ is Connes amenable
if and only if $G$ is virtually abelian.
\end{itemize}
\end{theorem}

\section{Operator amenability properties}

It is unfortunate that \falg\ is amenable so rarely.  The theory of operator
spaces, which developed into its modern form in the 1990s, gives a new
perspective to these problems.  Banach spaces (or even normed spaces)
have clear advantages over simple vector spaces for posing problems
in infinite dimensional settings.  The main reason is that there are useful
theorems regarding properties of bounded linear operators, while little can
be said about linear operators alone, in general.  The theory
of operator spaces and completely bounded maps offers a similar
refinement in terms of control, however still leaves us enough morphisms
for a usable theory.

\subsection{Operator spaces and completely bounded maps}
We use, as a standard reference for operator spaces, the book of 
Effros and Ruan~\cite{effrosrB}; all results in this section can be found there.  
We also recommend the book of Paulsen~\cite{paulsenB}.  

The axioms of Ruan, 
give a simple abstract definition of operator spaces.
Let $\fV$ be a $\Cee$-vector space and for each $n$ in $\En$, let
$\matn{\fV}$ denote the space of $n\cross n$ matrices with entries in $\fV$.
An {\it operator space structure} is a sequence of norms
$\bigl(\norm{\cdot}_n:\matn{\fV}\to\Ree^{\geq 0}\bigr)_{n\in\En}$ which satisfy
\[
\text{(OS1):  }\norm{\begin{bmatrix} U & 0 \\ 0 & V\end{bmatrix}}_{n+m}
=\max\{\norm{U}_n,\norm{V}_m\}, \qquad
\text{(OS2): } \norm{\alpha U\beta}_n\leq\norm{\alpha}\norm{U}_n\norm{\beta}
\]
where $U\in\matn{\fV},V\in\matm{\fV}$ and $\alpha,\beta\in\matn{\Cee}$,
where the latter space is normed as the space of linear operators on an 
$n$-dimensional Hilbert space.  We generally assume that $(\fV,\norm{\cdot}_1)$
is complete, and hence a Banach space.
We will call $\fV$, equipped with an operator
space structure, and {\it operator space}.
It is easily checked that any closed subspace of
$\fB(\fH)$ ($\fH$ a Hilbert space) is an operator space, where
we identify $\matn{\fB(\fH)}\cong\fB(\fH^n)$ in the usual manner.

A linear map between operator spaces $T:\fV\to\fW$ is called
{\it completely bounded} if its amplifications,
$[v_{ij}]\mapsto [Tv_{ij}]:\matn{\fV}\to\matn{\fW}$, are uniformly
bounded in $n$.  These are the natural morphisms of operator spaces.
We denote the space of such maps $\fC\fB(\fV,\fW)$.  It is itself an operator
space via the natural identification $\matn{\fC\fB(\fV,\fW)}
\cong\fC\fB(\fV,\matn{\fW})$, $[T_{ij}]\cong (v\mapsto [T_{ij}v])$.  
All bounded linear functionals are automatically 
completely bounded, thus $\fV^*$ is naturally an operator space.
In particular $\fsalg\cong\cstarg^*$ is an operator space and
$\falg$ inherits its operator space structure as a subspace.
We note that $\measg\cong\fC_0(G)^*$ is a {\em maximal} operator space,
as is the subspace \bloneg.  For a maximal operator space $\fM$, any bounded map
$T:\fM\to\fV$ is automatically completely bounded, whenever it is bounded.
In particular, operator space theory does not restrict the natural morphisms
emanating from \measg\ or from \bloneg.

A bilinear map $B:\fV\cross\fW\to\fZ$ on operator spaces is called
{\it jointly completely bounded} if the amplifications
$B^{(n,m)}:\matn{\fV}\cross\matm{\fW}\to M_{nm}(\fZ)$,
given by $B^{(n,m)}([v_{ij}],[w_{kl}])=[B(v_{ij},w_{kl})]$, are uniformly
bounded in pairs $(n,m)$.  The operator projective tensor product
$\fV\widehat{\otimes}\fW$ is the canonical object which linearizes
jointly completely bounded bilinear maps into completely bounded
bilinear maps.  A remarkable result, proved by Effros and Ruan
is a Grothendieck-type identification of the operator space projective
tensor product of the preduals of two von Nemann algebras:
if $\fM$ and $\fN$ are von Neumann algebras with preduals
$\fM_*$ and $\fN_*$, then $\fM_*\widehat{\otimes}\fN_*\cong
(\fM\overline{\otimes}\fN)_*$, the predual of the von Neumann tensor
product of $\fM$ and $\fN$.  It is not known to the author, a method of
proving this identity, without appeal to Tomita-Takeaski theory, in particular the
commutation formula $(\fM\overline{\otimes}\fN)'\cong\fM'\overline{\otimes}\fN'$.
For locally compact groups $G$ and $H$
this give rise to the identification 
\[
\falg\widehat{\otimes}\falh
\cong\fal{G\cross H}.  
\]
This contrasts with the result of Losert~\cite{losert},
which was proved, using techniques of subhomogeneous von Neumann algebras,
that $\falg\otimes^\gamma\falh\cong\fal{G\cross H}$, isomorphically, if an only if 
at least one of $G$ or $H$ is virtually abelian.  
Recall the role of virtually abelian groups
in Section \ref{sec:falgfsalba}.  These results suggest that
the recognition of the natural operator space structures on \falg\ and on \fsalg,
should play a role in gaining a more satisfactory amenability theory.
For example, using Loset's result, Forrest and 
Wood~\cite{forrestw} showed that all bounded maps on \falg\ are automatically
completely bounded exactly when $G$ is virtually abelian.

To this end, an algebra $\fA$, which is also an operator space, is called
a {\it completely contractive Banach algebra} if the multiplication map
$m:\fA\cross\fA\to\fA$ is jointly completely contractive, i.e.\ all amplifications $m^{(n,k)}$
are contractions.  An $\fA$-bimodule, $\fV$, which is also an operator space,
is a {\it completely contractive $\fA$-bimodule} provided the module
actions $\fA\cross\fV,\fV\cross\fA\to\fV$ are both jointly completely contractive,
which is the same as having that the each of the module-defining homomorphism and anti-
homomorphism from $\fA$ into $\fC\fB(\fV)$ are each completely contractive.
If $\fV$ is such a module, then $\fV^*$ is a dual completely contractive $\fA$-bimodule.
The entire theory of amenability, contractibility, weak amenability, biflatness
and biprojectivity (in the case of a function algebra, hyper-Tauberianness)
can be adapted to be rephrased with completely bounded derivations, 
operator projective tensor products, or module
maps (or local maps), and turned into the theories of {\it operator amenability},
{\it operator weak amenability}, etc.  In short, the adverb ``operator'' means
that we convert all results from the Banach space setting into the operator space
one.  Note that the maximal operator spaces structures on \bloneg\
and \measg\ allow the following important fact.

\begin{remark} Every Banach algebra theorem on \bloneg\
or on \measg\ is automatically a completely contractive
Banach algebra theorem.
\end{remark}

In other words, all of the results of Section \ref{sec:blonemeas}
are really operator amenability results.

\subsection{Fourier and Fourier-Stieltjes algebras as operator spaces}
Now, with operator spaces in hand, we can state the true analogues
of the theorems of Section \ref{sec:blonemeas}.  We let $G$ denote a locally
compact group.  The seminal
theorem is by Ruan~\cite{ruan}.

\begin{theorem}\label{theo:ruan}
\falg\ is operator amenable if and only if $G$ is amenable.
\end{theorem}

The success of this result has inspired constructions of operator space
structures on Herz--Fig\`{a}-Talamanca algebras $\mathrm{A}_p(G)$
to gain similar operator amenability results:  two quite different constructions
are obtained by Lambert et al~\cite{lambertnr} and Daws~\cite{daws1}.
The only systematic attempt, of which the author is aware, to 
generalize Theorems \ref{theo:johhel} and \ref{theo:ruan} to
a locally compact quantum group or Kac algebra setting is by Aristov~\cite{aristov04}.

Theorem \ref{theo:ruan} is the appropriate dual analogue to (i)$\Leftrightarrow$(ii) of
Theorem~\ref{theo:johhel}.  Since the dense inclusion $\falg\hookrightarrow\fC_0(G)$
is completely contractive, operator contractibility of \falg\ implies that of
$\fC_0(G)$.  Since each maximal ideal is co-dimension $1$, the arguments of 
\S 4.1 of \cite{rundeB} can be readily adapted to
show that $G$ is finite. The facts stated immediately after Theorem~\ref{theo:johhel},
along with duality considerations, lead to the following.

\begin{conjecture}\label{conj:opbiflat}
\falg\ is always operator biflat.
\end{conjecture}

We do not yet know the truth of this conjecture, but there is strong evidence in its 
favour.  We say that a closed subgroup $H$ in $G$ admits a {\it bounded
approximate indicator} in $G$ if there is a net $(v_\alpha)\subset\fsalg$ for which
\[
v_\alpha|_Hu\to u\text{ for }u\in\falh\quad\text{and}\quad
v_\alpha w\to 0\text{ for }w\in\ideal_{\falg}(H).
\]
(See explanation of notation in Section \ref{sec:commutative}.)
This concept was introduced in \cite{aristovrs}.  It was proved there that
bounded approximate indicators can always be built of positive definite
functions.  We observe that they may not be realized in general as elements
whose restriction to $H$ is the constant function $1$; see \cite{kaniuthl}.

We say a locally compact group $Q$ is a [QSIN] (quasi-small invariant 
neighbourhood) group if \blone{Q}\ admits a bounded
approximate identity $(e_\alpha)$ for which 
$\norm{e_\alpha\con\delta_s-\delta_s\con e_\alpha}_1$ $\to 0$
for each Dirac measure $\delta_s$ in $\meas{Q}$.

\begin{theorem}\label{theo:arisrunspr}
\begin{itemize}
\item[\rm (i)] If the diagonal subgroup $\Delta$ in $G\cross G$
admits a bounded approximate indicator, then \falg\ is operator biflat.
\item[\rm (ii)] If $G$ injects continuously into a [QSIN] group, then $\Delta$ admits
a bounded approximate indicator in $G\cross G$.
\item[\rm (ii)] If $G=\mathrm{SL}_3(\Ree)$, then $\Delta$ does not admit
a bounded approximate indicator in $G\cross G$.
\end{itemize}
\end{theorem}

This theorem was proved in \cite{aristovrs}, but
we note that Ruan and Xu~\cite{ruanx}
proved directly that the assumed condition in (ii), above, 
implies that \falg\ is operator biflat.
We do not know if $\fal{\mathrm{SL}_3(\Ree)}$ is operator biflat.

We observe that results (ii) and (i), with Leptin's theorem~\cite{leptin},
show directly the fact that for an amenable group $G$, \falg\ is operator
amenable.  It is interesting to note that a device in \cite{spronk}
allows us directly to obtain an operator bounded approximate diagonal
for \falg, for amenable $G$.  If $G$ is amenable then $G$ is itself
[QSIN] (see \cite{losertr84,stokke}), and hence, by (ii),
there is a bounded approximate indicator $(v_\alpha)$ for
$\Delta$ in $G\cross G$.  Let $(u_\beta)$ be a bounded
approximate identity for $\fal{G\cross G}$.  Then 
\[
(v_\alpha u_\beta)\subset\fal{G\cross G}\cong\falg\widehat{\otimes}\falg
\]
and a net $(w_\mu)$ can be extracted from this set
which allows limits to be taken first in $\alpha$, then in $\beta$.  This
net is the desired bounded approximate diagonal.

If Conjecture \ref{conj:opbiflat} were true, it would imply (i) of the following,
which is the dual analogue of Theorem~\ref{theo:bloneweakam}.

\begin{theorem}\label{theo:sprsam}
\begin{itemize}
\item[\rm (i)] \falg\ is always operator weakly amenable.
\item[\rm (ii)] \falg\ is always operator hyper-Tauberian.
\end{itemize}
\end{theorem}

Part (i) was proved in \cite{spronk}, and independently by Samei~\cite{samei05}.
While the proof of \cite{spronk} relied on 
the operator space analogue of Theorem~\ref{theo:spectral} (iii'), that
of \cite{samei05} was the genesis for the hyper-Tauberian condition of
\cite{samei06}.  It was there that (ii) was proved; moreover, this result was extended
to Herz--Fig\`{a}-Talamanca algebras $\mathrm{A}_p(G)$, using the
operator space structure of \cite{lambertnr}.  As mentioned
in Section \ref{sec:commutative}, (ii) implies (i).  
An important fact used in both of theses proofs is that subgroup $\Delta$ is
spectral for $\fal{G\cross G}$; see \cite{herz,takesakit}.
We note that it follows from Theorem
\ref{theo:sprsam}, and the fact that bounded linear functionals are automatically
completely bounded, that \falg\ admits no non-zero point derivations.
 
We have the following dual analogue of Theorem~\ref{theo:blonebiproj}
which is due to Aristov~\cite{aristov} and, independently, Wood~\cite{wood}.
 
 \begin{theorem}\label{theo:opbiproj}
 \falg\ is operator biprojective if and only if $G$ is discrete.
 \end{theorem}

Daws~\cite{daws2} has recently provided evidence which suggests
that Theorems \ref{theo:blonebiproj} and \ref{theo:opbiproj}
are really results of compact Kac algebras, and may not be extended
to general compact quantum groups.

The success of finding the dual analogue of each amenability result about 
\bloneg\ in Section \ref{sec:blonemeas} for \falg\ (except, possibly
operator biflatness) lead us to think that Theorem~\ref{theo:daghhe}
should admit a dual analogue for \fsalg\ in the operator space setting.
Thus (ii) and (iii) below, were a surprise.  The operator amenability
constant $C_\fA^{op}$, is defined for a completely contractive Banach algebra
$\fA$ analoguously to the amenability constant of Section \ref{ssec:amenability}.

\begin{theorem}\label{theo:fell}
\begin{itemize}
\item[\rm (i)] If \fsalg\ is operator amenable with operator amenability 
constant $C_{\fsalg}^{op}<5$, then $G$ is compact.
\item[\rm (ii)] There is a class of groups, the ``Fell'' groups, $G$, for which
each \fsalg\ is operator amenable with $C_{\fsalg}^{op}=5$.
\item[\rm (iii)] For each Fell group $G$, \fsalg\ is weakly amenable.
\end{itemize}
\end{theorem}

Part (i) was proved in \cite{rundes03}, while parts (ii) and (iii) are from 
\cite{rundes07}.  For a Fell group $G$ -- an example of which is the 
semi-direct product $G_p=\mathbb{Z}_p\ltimes\mathbb{Q}_p$ for a prime $p$ --
we have the decomposition
$\fsalg=\ffalg\oplus_{\ell^1}\falg$.  This is a ``semi-direct product''
of two operator amenable Banach algebras.  The lack of the adverb
``operator'' in (iii) is not a typographical error.  Since a Fell group $G$
is totally disconnected, and so too is its almost periodic compactification
-- in the case $G_p$, above, it is $\mathbb{Z}_p$ -- each algebra
$\falg$ and $\ffalg\cong\fal{G^{ap}}$ is weakly amenable by
Theorem~\ref{theo:runfor} (ii).  

In particular we notice that for a Fell group, \fsalg\ has no non-zero point 
derivations.
There are other known examples of this absence of point derivations.
If $G$ is a connected semi-simple Lie group, it was shown by
Cowling~\cite{cowling} that there is a finite family of closed
normal subgroups $N_1,\cdots,N_n$, with quotient maps
$q_j:G\to G/N_j$, for which $\fsalg=\ell^1\text{-}\bigoplus_{j=1}^n
\fsoal{G/N_j}\comp q_j$.  Moreover, Cowling applied his
generalized Kunze-Stein phenomenon~\cite{cowlingKS},
to show that each $\fsoal{G/N_j}/\fal{G/N_j}$ is a radical Banach algebra.
Since each $\fal{G/N_j}$ admits no non-zero point derivations (Theorem
\ref{theo:sprsam}), neither can $\fsoal{G/N_j}$.  These examples
stand in marked contrast to (i)$\Leftrightarrow$(ii) of Theorem~\ref{theo:daghhe}.
(The author is grateful to M.\ Ghandehari for pointing out this example.)

It is further interesting to consider amenability properties of the Rajchman
algebras $\fsoalg$.  In contrast to the examples above, if
$G$ is abelian and non compact, Ghandehari~\cite{ghandehari},
has shown that point derivations exist on $\fsoalg$.


\begin{question}
\begin{itemize}
\item[\rm (i)] When is \fsalg\ operator amenable? Operator weakly amenable?
\item[\rm (ii)] What possible values are there for $C^{op}_{\fsalg}$?
\item[\rm (iii)] When is \fsoalg\ operator amenable? Operator weakly amenable?
\item[\rm (iv)]  When do \fsalg, or \fsoalg, admit non-zero point derivations?
\end{itemize}
\end{question}

In partial answer to (ii),
it is suggested in \cite{ghandeharihs} that  $C^{op}_{\fsalg}$ takes
only values $4n+1$.  However no example is known to us of examples taking values
other than $1$ or $5$.  It is reasonable to expect that $\fsoalg$ is operator
weakly amenable only if $\fsoalg=\falg$.

Even when the dual structure is taken into account, the situation for \fsalg\
is still not clear.

\begin{theorem}\label{theo:fsalgconnesam}
\begin{itemize}
\item[\rm (i)] \fsralg\ is operator Connes amenable if and only if $G$ is amenable.
\item[\rm (ii)] If $G=\mathrm{F}_2$, the free group on two generators, then
\fsalg\ is operator Connes amenable.
\end{itemize}
\end{theorem}

This theorem is proved in \cite{rundes03}.  The proof of (i) is easy:
$\fsralg$ can be operator Connes amenable only if it is contains the
constant function $1$; conversely if $G$ is amenable, then $\falg$ is
operator amenable and weak*-dense in \fsralg.  The proof of (ii)
invokes the fact, due to Choi~\cite{choi}, that for $G=\mathrm{F_2}$,
$\cstarg$ is residually finite dimensional. This implies that
\ffalg, which is operator amenable, is weak* dense in \fsalg.
We have, presently, 
no natural conjectures as to when \fsalg\ is operator Connes amenable,
in general.   It is not clear that \fsalg\ ought to admit an operator normal virtual 
diagonal in such a case; this condition is not automatic \cite{runde06a}.

\section{Applications}

\subsection{Ideals of Fourier algebras with bounded approximate identities}
\label{sec:ideals}
It is a natural question, given a Banach algebra, to ask if it is possible
to classify its ideals.  For an algebra such \falg, even for $G=\mathbb{T}$,
this is a difficult question, since it is known that there are non-spectral subsets
(see definition in Section \ref{sec:commutative}).  Hence one may wish
to classify a nice subclass of ideals, say, those admitting bounded approximate
identities. 
For abelian groups this was achieved by Liu et al~\cite{liuvrw},
using methods of invariant means; this pointed towards the need for some
form of averaging in \falg.  Since \falg\ is amenable only for virtually abelian
$G$, this points to the necessity to import operator space techniques so
we can access arbitrary amenable $G$.   

Let $\Omega(G)$ denote the smallest Boolean ring of subsets
containing cosets of subgroups, and $\Omega_c(G)$ the
set of closed elements in $\Omega(G)$.  The following is the main
result of \cite{forrestkls}.

\begin{theorem}\label{theo:idealbai}
Let $G$ be amenable.  Then every closed ideal of \falg, which has a bounded
approximate identity, is of the form $\ideal_{\falg}(E)$ where
$E\in\Omega_c(G)$.  
\end{theorem}

The proof is based on the fact that an ideal $\fI$, in an (operator)
amenable Banach algebra $\fA$, has a bounded approximate identity
if [and only if] it is weakly (completely) complemented. i.e.\ there
is a (completely) bounded projection $P:\fA^*\to\fI^\perp$, the annihilator
of $\fI$ in $\fA^*$ (see \cite{helemskii,rundeB}).  This follows from the fact
that this projection can be ``averaged'' to a (completely) bounded $\fA$-bimodule
projection.  Now if $H$ is a closed subgroup of an amenable group $G$, then
$\vnh$ is an injective von Neumann algebra, from which it can be deduced
that there is a completely bounded projection $P:\vng\to\ideal_{\falg}(H)^\perp
\cong\vnh$.  Since, by Theorem~\ref{theo:ruan}, $\falg$ is operator amenable,
$P$ can be averaged to invariance
(see \cite{wood00,forrests}), hence giving the fact
that $\ideal_{\falg}(H)$ admits a bounded approximate identity. Then, careful
arguments, using the structure of elements of $\Omega_c(G)$, shows that for any
$E$ in $\Omega_c(G)$, $\ideal_{\falg}(E)$ has a bounded approximate identity.

To see the converse, consider $\fI$ a closed ideal in \falg\
with bounded approximate identity $(u_\alpha)$, and let 
$E=\{s\in G:u(s)=0\text{ for }u\in\fI\}$.  Then $(u_\alpha)$
has a weak* (pointwise) cluster point in $\fsal{G_d}$ ($G_d$
is $G$ with discrete topology) which is necessarily the indicator function
$1_{G\setminus E}$.  By Host~\cite{host}, $E\in\Omega(G)$, and is closed,
so $E\in\Omega_c(G)$.  We, moreover, show that elements of
$\Omega_c(G)$ are spectral sets, hence $\fI=\ideal_{\falg}(E)$.

We remark that in the setting of Fourier algebras important results leading to
Theorem \ref{theo:idealbai} were obtained by Forrest~\cite{forrest90,forrest92}.
Operator space techniques were applied to this problem for the first time
by Ruan and Xu~\cite{ruanx} and Wood~\cite{wood00}.

In the case that $G$ is discrete, Wood~\cite{wood03} uses operator
biprojectivity (Theorem~\ref{theo:opbiproj}) to obtain the following.

\begin{theorem}
If $G$ is discrete, then a closed ideal $\fI$ in $\falg$
is completely complemented if and only if there is a completely bounded
$\falg$-bimodule map $P:\falg\to\fI$.
\end{theorem}

In the case that $G$ is amenable it follows that $\fI=\ideal_{\falg}(E)$
for some $E$ in $\Omega(G)$.  However, for non-commutative free groups $F$ there
are certain free sets $E$, due to Leinert~\cite{leinert}, which are not elements of
$\Omega(F)$, but for which $\ideal_{\fal{F}}(E)$ is completely complemented.

The desire to understand ideals in $\falg$ for non-amenable $G$ has led
us to consider the algebra $\fcbal{G}$, which is the closure of
\falg\ in the algebra of completely bounded multipliers $\multcb\falg$;
the latter space was introduced by De Canniere and Haagerup~\cite{decanniereh}
and shown to be a dual space.  For ``weakly amenable'' $G$,
as defined in \cite{decanniereh}
(this should not be confused with weakly amenable algebras),
we always have that $\fcbal{G}$ admits a bounded approximate
identity (see \cite{forrest05}).  Thus, $\fcbal{G}$ represents faithfully
on many canonical modules, such as $\falg$ and \vng, whereas \falg\ does not.
It is shown in \cite{decanniereh} that $\mathrm{SL}_2(\Ree)$ is weakly amenable
and that $\mathrm{F_2}$, being a finite index subgroup of a lattice in 
$\mathrm{SL}_2(\Ree)$, also enjoys this property.
The following result is from \cite{forrestrs} and uses techniques similar
to those in Theorem~\ref{theo:fsalgconnesam} (ii).  

\begin{theorem}
\begin{itemize}
\item[\rm (i)] $\multcb\fal{\mathrm{F_2}}$ is operator Connes amenable.
\item[\rm (iii)] $\fcbal{\mathrm{F_2}}$ is operator amenable.
\item[\rm (iii)] Every weakly complemented closed ideal in $\fcbal{\mathrm{F_2}}$
has a bounded approximate identity.
\end{itemize}
\end{theorem}

\subsection{Fourier algebras of homogeneous spaces}
Let $G$ be a locally compact group and $K$ a compact subgroup.
The Fourier algebras of the homogeneous spaces 
$\falgk$ were defined by Forrest~\cite{forrest98}:
$\falgk=\{u\in\falg:k\mult u=u\text{ for }k\in K\}$, where
$k\mult u(s)=u(sk)$.  It was shown, in that article, that these
spaces have many of the accoutrements of operator amenable algebras:
for example, they admit bounded approximate identities when
$G$ is amenable.  Moreover, the analogue of the
Grothendieck-Effros-Ruan tensor product formula was shown in
\cite{forrestss1}:  $\falgk\widehat{\otimes}\falgk\cong
\fal{G/K\times G/K}$ (where, of course $G/K\cross G/K$ is the homogeneous
space $(G\cross G)/(K\cross K)$).  Hence the following contrast to 
Theorem \ref{theo:ruan}, proved in \cite{forrestss1}, was a surprise.

\begin{theorem}\label{theo:forrestss}
If $G$ is a compact semi-simple Lie group, then
$\fal{G\cross G/\Delta}$ is not operator weakly amenable.
\end{theorem}

As an interesting consequence, we could show that
for $G$, as in the theorem, $(G\cross G\cross G\cross G)\Delta\cross\Delta$
is not spectral for $\fal{G\cross G\cross G\cross G}$.  
We note that the natural identification $s\mapsto (s,e)\Delta$
gives rise to an isomorphic identification $\fal{G\cross G/\Delta}\cong\falg$
only when $G$ is virtually abelian.  Otherwise $\fal{G\cross G/\Delta}$
identifies with a subalgebra $\fdelg$ of \falg, which is very similar to 
the algebra $\fgamg$ of Johnson~\cite{johnson94}.  Hence Theorem \ref{theo:forrestss}
is really an analogue of Theorem \ref{theo:falgnwa}.

We say that $G$ is [MAP]$_K$
if the almost periodic functions on $G$ separate the points on $K$.
In \cite{forrestss1} we obtained the following results for locally compact $G$.

\begin{theorem}
\begin{itemize}
\item[\rm (i)] If $G_e$ is abelian, then $\falgk$ is hyper-Tauberian, and
if $G$ is [MAP]$_K$, then
$(K\cross K)\Delta$ is a set of local synthesis for $\fal{G\cross G}$.
\item[\rm (ii)] If $G$ is amenable and [MAP]$_K$, then $\falgk$ is operator amenable
if and only if $(K\cross K)\Delta\in\Omega(G\cross G)$.
\end{itemize}
\end{theorem}

If $G$ contains a compact non-abelian connected subgroup $K$, then
there is a compact subgroup $K^*$ of $G\cross G$ for which
$\fal{G\cross G/K^*}$ is not operator weakly amenable; a partial converse to (i).
The only cases which we know of where $(K\cross K)\Delta\in\Omega(G\cross G)$
are when $K$ has a subgroup of finite index which is normal in $G$.

In the follow-up paper \cite{forrestss2} we gained an improvement of Theorem
\ref{theo:forrestss}, but with less illustrative methods.

\begin{theorem}
Let $G$ be a compact group.  Then the following are equivalent:
\begin{itemize}
\item[\rm(i)]  $\fdelg$ is operator weakly amenable;
{\rm (i')} $\fdelg$ is weakly amenable;
\item[\rm(ii)] $\fdelg$ is operator hyper-Tauberian;
{\rm (ii')} $\fdelg$ is hyper-Tauberian;
\item[\rm (iii)] $G_e$ is abelian.
\end{itemize}
Moreover, the following are equivalent:
\begin{itemize}
\item[\rm(a)] $\fdelg$ is operator amenable; 
{\rm (a')} $\fdelg$ is amenable;
\item[\rm (b)] $G$ is virtually abelian;
{\rm (b')} $\fdelg=\falg$.
\end{itemize}
\end{theorem}

What is interesting about this result, is that it suggests that $\fdelg$ contains,
in its Banach space structure, critical information about the operator space
structure of \falg.

\subsection{Homomorphisms on Fourier algebras}
Let $G$ and $H$ be locally compact groups.
Cosets of $H$ can be characterized at those subsets $C$
which are closed under the ternary operation $(r,s,t)\mapsto rs^{-1}t$.
A map $\alpha:H\to G$ is called {\it affine} if $\alpha(rs^{-1}t)
=\alpha(r)\alpha(s)^{-1}\alpha(t)$ for $r,s,t$ in $C$.  
We say $\alpha:Y\to G$ is {\it piecewise affine}
if there is a partition $Y=\bigcup_{j=i}^nY_n$ of $Y$, 
where each $Y_j\in\Omega(H)$ and for each $j$
a coset $C_j\supset Y_j$ and an affine $\alpha_j:Y_j\to G$ for which
$\alpha|_{Y_j}=\alpha_j|_{Y_j}$.  We will call $\alpha:Y\subset H\to G$
a {\it continuous piecewise affine map} if each $Y_j\in\Omega_o(G)$
(the smallest ring of subset generated by open cosets) and $\alpha$ is continuous.
If $u:G\to\Cee$ is a function, we can define for piecewise affine
$\alpha$, as above, the map $\Phi_\alpha u:H\to\Cee$ by
$\Phi_\alpha u(s)=u(\alpha(s))$ if $s\in Y$, and $\Phi_\alpha u(s)=0$
if $s\in H\setminus Y$.

The following result was proved, for abelian groups by Cohen~\cite{cohen},
generalized to the case that $G$ is virtually abelian by Host~\cite{host}, and
then to the case that $G$ is discrete and amenable by Ilie~\cite{ilie}.
The definitive form, as presented here, is in \cite{ilies}.

\begin{theorem}\label{theo:iliespr}
\begin{itemize}
\item[\rm (i)] If $\alpha:Y\subset H\to G$ is continuous piecewise affine,
then $\Phi_\alpha:\falg\to\fsalh$ is a completely bounded homomorphism.
\item[\rm (ii)] If $G$ is amenable, then every completely bounded
homomorphism $\Phi:\falg\to\fsalg$ is of the form $\Phi=\Phi_\alpha$, as in {\rm (i)}, above.
\item[\rm (iii)] The homomorphism $\Phi_\alpha$ maps \falg\ into \falh\
if and only if $\alpha$ is proper, i.e.\ $\alpha^{-1}(K)$ is compact
for each compact $K\subset H$.
\end{itemize}
\end{theorem}

As an obvious corollary, we have that if $H$ is connected, only affine
$\alpha$ are allowed as ``symbols'' for completely bounded homomorphisms
from \falg\ ($G$ amenable) to \fsalg.  This theorem seems nearly sharp
in the following respects.  First, the map $u\mapsto\check{u}:\falg\to\falg$
($\check{u}(s)=s(s^{-1})$) is completely bounded only if $G$ is virtually
abelian; compare to Theorem \ref{theo:runfor} (ii) and comments in the 
paragraph below.   Second, if $E$ is a free set
in a non-commutative, hence non-amenable, free group $F$, then 
$u\mapsto 1_F u:\fal{F}\to\fal{F}$
is completely bounded (see Leinert~\cite{leinert}) but not implemented
by a completely bounded symbol (in particular $E\not\in\Omega(F)$).

Let us illustrate the role of Theorem \ref{theo:ruan} in the proof of (ii).
If $\Phi$ is completely bounded, then 
\[
\Phi\otimes\id_{\falg}:\falg\widehat{\otimes}\falg\cong\fal{G\cross G}
\to\fsal{H\times G}.
\]
If $(w_\mu)$ is an operator bounded approximate diagonal for \falg\
of the type whose construcion is outlined after Theorem \ref{theo:arisrunspr},
then the net $(\Phi(w_\mu))$, considered as a bounded net in
$\fsal{H_d\otimes G_d}$ (discretized groups) must have a cluster point
which is the indicator function $1_Y$.  Hence $Y\in\Omega(H\cross G)$,
by Host~\cite{host}.  It can then be checked that $Y$ is the graph of a piecewise
affine function $\alpha$, which, moreover, must be continuous.

We define $\alpha:Y\subset H\to G$ to be {\it mixed piecewise affine}
if $Y$ partitions as above into elements $Y_1,\dots,Y_n$ of $\Omega(H)$,
and for $j=1,\dots,n$ there is a coset $C_j\subset Y_j$ and either an
affine $\alpha_j:C_j\to G$, or an anti-affine $\alpha_j$ (i.e.\ $\alpha_j(rs^{-1}t)
=\alpha_j(t)^{-1}\alpha_j(s)\alpha_j(s)^{-1}$), such that 
$\alpha|_{Y_j}=\alpha_j|_{Y_j}$.  If $\Phi_\alpha$ is define as above,
it is clear that $\Phi_\alpha:\falg\to\fsalg$ is a bounded homomorphism.

\begin{conjecture}
If $G$ is amenable and $\Phi:\falg\to\fsalh$ every bounded homomorphism,
is $\Phi=\Phi_\alpha$ for some continuous mixed piecewise affine $\alpha:Y\subset G
\to H$?
\end{conjecture}

It is an interesting open question to characterize completely bounded, never mind
bounded, homomorphisms on Fourier algebras of non-amenable groups.
As suggested by the proof Theorem \ref{theo:iliespr}, illustrated
above, the completely bounded question is closely
linked to our understanding of the idempotents in $\multcb\fal{H\cross G}$
(see section \ref{sec:ideals} above for notation).

For group algebras, only the contractive homomorphisms
$\Psi:\bloneg\to\meas{H}$ are well understood; see Greenleaf~\cite{greenleaf}.
It would be interesting to see if, for amenable $G$, bounded approximate
diagonal methods shed light on the understanding of bounded such
homomorphisms.  If these structures can be understood, then there is
hope for understanding the completely bounded homomorphisms on
predual algebras of locally compact quantum groups.

\bigskip\noindent 
{\bf Acknowledgement.}

The author's research is supported by NSERC, under grant No.\ 312515-05.

This paper is based on a lecture delivered at the 19$^\mathrm{th}$ 
International Conference on Banach Algebras held at B\c{e}dlewo, July 14--24, 
2009. The support for the meeting by the Polish Academy of Sciences, the 
European Science Foundation under the ESF-EMS-ERCOM partnership, and the 
Faculty of Mathematics and Computer Science of the Adam Mickiewicz University 
at Pozna\'n is gratefully acknowledged.

The author is grateful to E.\ Samei for reading the manuscript and
correcting some errors, and to the referee for correcting many typos
and making suggestions to improve the exposition.

\end{document}